\documentclass[12pt]{amsart}
\usepackage{amsmath,amsfonts,amssymb,amsxtra,latexsym,amscd,enumerate,amsthm}

\usepackage{latexsym}
\usepackage{epsfig}

\def\q{\frac{1}{2}}

\def\e{\varepsilon}

\def\DD{\Delta}

\def\o{\omega}
\newcommand\D{{\mathcal D}}
\def\dh{{\mathcal{D}H}}
\def\dl{{\mathcal{D}L}}
\def\dx{{\mathcal{D}X}}
\def\dy{{\mathcal{D}Y}}

\def\S{\mathbb{S}}
\def\R{\mathbb{R}}
\def\C{\mathbb{C}}
\def\Z{\mathbb{Z}}
\def\N{\mathbb{N}}
\def\H{{\mathcal H}}
\def\les{\lesssim}

\def\beq{\begin{equation}}
\def\eeq{\end{equation}}

\def\beq{\begin{equation}}
\def\eeq{\end{equation}}

\newtheorem{d1}{Definition}
\newtheorem{t1}{Theorem}
\newtheorem{l1}{Lemma}
\newtheorem{p1}{Proposition}
\newtheorem{c1}{Corollary}

\begin{document}

\title[Large data local solutions for the D-NLS equation]{Large data
  local solutions for the derivative NLS equation}

\author{Ioan Bejenaru}
\address{ Department of Mathematics \\
  University of California, Los Angeles}

\author{ Daniel Tataru}
\address {Department of Mathematics \\
  University of California, Berkeley}

\thanks{ The first author was
  partially supported by MSRI for Fall 2005 \\ The second author was
  partially supported by NSF grants DMS0354539 and DMS 0301122 and
  also by MSRI for Fall 2005}

\begin{abstract}
  We consider the Derivative NLS equation with general quadratic
  nonlinearities. In \cite{be2} the first author has proved a sharp
  small data local well-posedness result in Sobolev spaces with a
  decay structure at infinity in dimension $n = 2$. Here we prove a
  similar result for large initial data in all dimensions $n \geq 2$.

\end{abstract}

\maketitle

\section{Introduction}

The general Cauchy problem for the semilinear Schr\"odinger equation  
has the form
\begin{equation} \label{E}
\begin{cases}
\begin{aligned} 
  iu_{t}-\Delta u &= P(u,\bar{u},\nabla u, \nabla \bar{u}), 
\   t\in \mathbb{R}, x\in \mathbb{R}^{n} \\  u(x,0) &=u_{0}(x)
\end{aligned}
\end{cases}
\end{equation}
where $u: \mathbb{R}^{n} \times \mathbb{R} \rightarrow \mathbb{C}^m$ and
$P: \mathbb{C}^{2n+2} \rightarrow \mathbb{C}^m$.  
Assuming that $P$ is a polynomial containing terms of order at least
$\kappa \geq 2$ and higher, one is interested in studying the local
well-posedness for this evolution in a suitable Sobolev space. 

In the simpler case when $P$ does not depend on $\nabla u, \nabla
\bar{u}$ this problem is rather well understood, and the Strichartz
estimates for the linear Schr\"odinger equation play a proeminent
role. 

Here we are interested in nonlinearities which contain derivatives.
This problem was considered in full generality in the work of
Kenig-Ponce-Vega~\cite{k2}. Due to the need to regain one derivative in
multilinear estimates, they use in an essential fashion the
local smoothing estimates for the linear Schr\"odinger equation.
Their results make it clear that one needs to differentiate two cases.

If $\kappa \geq 3$ then they prove local well-posedness for initial
data in a Sobolev space $H^N$ with $N$ sufficiently large. However,
if quadratic nonlinearities are present, i.e. $\kappa = 2$, then the 
local well-posedness space also incorporates decay at infinity, namely
\[
H^{N,N} = \{ x^N u \in L^2; \ D^N u \in L^2\}
\]  
with $N$ large enough.  The need for decay is motivated by work of 
Mizohata \cite{m}.  He proves that a necessary condition   for the $L^{2}$ well-posedness
of the problem:
\begin{equation}
\begin{cases}
\begin{aligned} 
iu_{t}-\Delta u &= b_{1}(x) \nabla u , \   t\in \mathbb{R}, x\in
\mathbb{R}^{n} 
\\u(x,0) &=u_{0}(x)
\end{aligned}
\end{cases}
\end{equation}
is the uniform bound
\beq 
\label{co1}
\sup_{x \in \R^{n}, \omega \in \mathbb{S}^{n-1}, R > 0} \left|Re
  \int_{0}^{R} b_{1}(x+r \omega) \cdot \omega dr\right | < \infty \eeq
The idea behind this condition is that $Re \ b_{1}$ contributes to
exponential growth of the solution along the Hamilton flow of the
linear Schr\"odinger operator. 

If $\kappa \geq 3$ then naively one needs such bounds for expresions
which are quadratic in $u$. If $u \in H^{N}$ then $u^2 \in W^{N,1}$
and the above integrability can be gained. However, if $\kappa = 2$ 
then one would want similar bounds for linear expressions in $u$;
but this cannot follow from square integrability, therefore one compensates
by adding decay at infinity.

In this work we consider the case of quadratic nonlinearities, precisely
the problem
\begin{equation} \label{Eg}
\begin{cases}
\begin{aligned} 
  iu_{t}-\Delta u &= B((u,\bar u),(\nabla u,\nabla \bar u)),  \\  u(x,0) &=u_{0}(x)
\end{aligned}
\end{cases}
\end{equation}
Here $B$ is a generic bilinear form which contains one differentiated
and one undifferentiated factor, and also may have complex conjugates. Our
results easily transfer via differentiation to the similar problem
with two derivatives in the nonlinearity,
\begin{equation} \label{Egg}
\begin{cases}
\begin{aligned} 
  iu_{t}-\Delta u &= B((\nabla u,\nabla \bar u), (\nabla u,\nabla \bar
  u)), \\ u(x,0) &=u_{0}(x)
\end{aligned}
\end{cases}
\end{equation}

The initial data $u_0$ is assumed to be locally  in Sobolev spaces $H^s$ 
but with some additional decay at infinity. In what follows $\lambda$
is a dyadic index. For a function $u$ we consider a Littlewood-Paley 
decomposition in frequency
\[
u = \sum_{\lambda \geq 1} u_\lambda, \qquad u_\lambda = S_\lambda u
\]
where all the frequencies smaller than $1$ are included in $u_1$.
Then following \cite{be2} we define the spaces $\dh^s$ by
\[
\| u\|_{\dh^s}^2 = \sum_{\lambda  \geq 1} \lambda^{2s} \|u_\lambda\|_{\dl^2_\lambda}^2
\]
where the dyadic norms $\dl^2_\lambda$ are defined in a manner somewhat 
similar to \eqref{co1}, namely
\[
  \|v\|_{\dl^2_\lambda} = \sup_{x_0 \in \R^n} \sup_{\omega \in \S^{n-1}}
\sum_{k \in \N} \| 1_{\{|\lambda^{-1} (x-x_0)-k\omega|<1\}} u \|_{L^2}.
\]
This definition is consistent with the speed of propagation properties
for the linear Schr\"odinger equation. Waves with frequency $\lambda$ have
speed $\lambda$ therefore move about $\lambda$ within a unit time interval.
Hence the linear Schr\"odinger equation is well-posed\footnote{The
  operator norm of the evolution in $\dh^s$ will grow polynomially in
time though.} in $\dh^s$.

If one considers the low regularity well-posedness for \eqref{Eg},
respectively \eqref{Egg} in $\dh^s$ then a natural threshold is given
by scaling, namely $s_c = \frac{n}2-1$ for \eqref{Eg} respectively
$s_c = \frac{n}2$ for \eqref{Egg}. However, it turns out that the
obstruction identified in Mizohata's work is much stronger and leads
to additional restrictions.  In \cite{ch1} Chihara obtains some better
results on this problem, lowering the threshold for \eqref{Egg} to $s=
\frac{n}2+4$. However, this is still far from optimal.  Indeed, a
sharp small data result is obtained by the first author in a recent
paper:

\begin{t1} [Bejenaru~\cite{be2}]
If $n=2$,  the equations \eqref{Eg}, respectively \eqref{Egg} are locally
  well-posed for initial data which is small in
  $\dh^{s}$ for all $s > s_c+1$.
  \end{t1}

  Although the above theorem was proved in dimension two, the same
  result can be derived in all dimensions $n \geq 2$.  In addition, if
  some limited spherical symmetry is imposed on the data then the
  above exponents can be relaxed up to scaling, see \cite{be1}.

The goal of this paper is to obtain a local well-posedness result for
the same initial data space as in \cite{be2}, but for large initial
data.  Our main result is

\begin{t1}
  Let $n \geq 2$ and $s > s_c+1$. Then the equations \eqref{Eg},
  respectively \eqref{Egg} are locally well-posed for initial data in
  $\dh^s$.  More precisely, there is $C > 0$ so that for each initial
data  $u_0 \in \dh^s$ there is an unique solution 
\[
u \in C(0,T;\dh^s), \qquad T = e^{-C\|u_0\|_{\dh^s}}
\]
In addition, the solution has a Lipschitz dependence on the initial data.
\label{maint}\end{t1}

We believe that this result is sharp for generic nonlinearities.
However, there are special cases when one is able to obtain stronger
results, see for instance \cite{co}, \cite{co2}, \cite{gr} ,\cite{ha}, \cite{st}.
 These results consider nonlinearities of the type $B((u,\bar{u}), \nabla{\bar{u}})$ and $B(\nabla{\bar{u}}, \nabla{\bar{u}})$ 
 in various dimensions and show that the regularity threshold for the initial data can be lowered all the way to the scaling.

The step from small to large data is entirely nontrivial.  The
difficulty is related to the infinite speed of propagation for the
linear Schr\"odinger equation. Precisely, a large low frequency
component of the solution produces an exponentially large perturbation
in the high frequency flow even for an arbitrarily short time.  Thus
one needs to add the low frequency part of the data to the linear
equation and only then do a perturbative analysis for the high
frequencies.  A somewhat similar analysis has been carried out before
for related problems, see \cite{k2}, \cite{ch1}; however, in both
works the full initial data becomes part of the nonlinearity, leading
to considerable technical difficulties.  What also differentiates the
present work is that the perturbed linear equation is very close to
Mizohata's necessary condition.

Another interesting feature of this work is the choice of the function
spaces for the perturbative analysis.  It has been known for some time
that the $X^{s,b}$ spaces are not good enough in order to study the
local theory for derivative NLS equation.  This is due to the need to
regain one derivative in the bilinear estimates, which leads to
logarithmic losses.  In \cite{be2} the first author introduced a
refinement of the $X^{s,b}$ spaces which removes this difficulty. In
this article we provide an alternate modification of the $X^{s,b}$
spaces which seems better suited for the study of variable coefficient
equations. This is based on a wave packet type decomposition
of solutions on the uncertainty principle scale.

\section {Scaling and the perturbative argument}

Differentiating once the equation \eqref{Egg} we obtain a equation of
type \eqref{Eg}. Thus in what follows we restrict our analysis to
\eqref{Eg}.  We denote
\[
M=\|u_{0}\|_{\mathcal{D}H^{s}}
\]
Our equation is invariant with respect to the scaling

\beq \label{f1} u^{\e}(x,t)=\e u (\e x,\e^{2}t), \ \ u^{\e}_{0}(x)=\e
u_{0}(\e x) \eeq It is natural to seek to decrease the size of the
initial data is by rescaling with a sufficiently small parameter $\e$.
However, there is a difficulty arising from the fact that we are using
inhomogeneous Sobolev spaces. This is why the result in \cite{be2}
does not yield a result for the large initial data.

We split the rescaled initial data $u_0^\e$ into low and high
frequencies,
\[
u^{\e}_{0} = u^{\e}_{0, \leq 1} + u^{\e}_{0, > 1}
\]
It is not too difficult to estimate their size:

\begin{p1} \label{low1} Assume that $s > \frac{n}2$. Then the
  components $u^{\e}_{0, \leq 1}$ and $u^{\e}_{0, > 1}$ of $u^\e_0$
  satisfy the pointwise bounds
  \begin{equation}
    \|u^{\e}_{0, \leq 1}\|_{L^\infty}+ \|u^{\e}_{0, > 1}\|_{L^\infty}
    \lesssim \e M
    \label{i2}\end{equation}
  and the $L^2$ estimates\footnote{The constant in the second bound
    needs to be adjusted to $\e |\ln\e|$ if $s = \frac{n}2 +1$. }
  \begin{equation}
    \| u^{\e}_{0, \leq 1}\|_{\dl^2} \lesssim M, \qquad \| \nabla
    u^{\e}_{0, \leq 1}\|_{\dl^2} \lesssim \max\{\e,\e^{s-\frac{n}2}\} M,
    \label{ulow}\end{equation}
  respectively
  \begin{equation}
    \| u^{\e}_{0, > 1}\|_{\dh^s} \lesssim \e^{s-\frac{n}2} M.
    \label{uhigh} \end{equation}
\end{p1}

The low frequency component is large but has the redeeming feature
that it does not change much on the unit time scale; hence we freeze
it in time modulo small errors. On the other hand the high frequency
component is small, therefore we can treat it perturbatively on a unit
time interval. To write the equation for the function
\[
v^{\e}=u^{\e}-u_{0,\leq 1}^{\e}
\]
we decompose the quadratic nonlinearity $B$ depending on whether the
gradient factor is complex conjugate or not,
\[
B((u,\bar u),(\nabla u,\nabla \bar u)) = B_0((u,\bar u),\nabla u) +
B_1((u,\bar u),\nabla \bar u)
\]
Then $v^\e$ satisfies

\begin{equation} \label{e1} \left\{
    \begin{aligned} 
      & iv_{t}-\Delta v - A v = B((v,\bar v), (\nabla{v},\nabla \bar
      v))+ A_1 v + N(u^{\e}_{0, \leq 1}) \\ &v(0) =u^{\e}_{0, > 1}
    \end{aligned}
  \right.
\end{equation}
where
\[
A v = B_0((u_{0,\leq 1}^{\e}, \bar u_{0,\leq 1}^{\e}), \nabla{v})
\]
is a linear term which is included in the principal part,
\[
A_1 v = B((v,\bar v), (\nabla u^{\e}_{0, \leq 1},\nabla \bar
u^{\e}_{0, \leq 1})) + B_1((u^{\e}_{0, \leq 1},\bar u^{\e}_{0, \leq
  1}),\nabla \bar v)
\]
is a linear term which can be treated perturbatively and
\[
N(u^{\e}_{0, \leq 1}) = B((u^{\e}_{0, \leq 1}, \bar u^{\e}_{0, \leq
  1}) , (\nabla u^{\e}_{0, \leq 1}, \bar \nabla u^{\e}_{0, \leq 1})) +
\DD u^{\e}_{0, \leq 1}
\]
represents the time independent contribution of the low frequency part
of the data.

To solve this we need some Banach spaces $X^s$, $\dx^s$ for the
solution $v$ respectively $Y^s$, $\dy^s$ for the inhomogeneous term in
the equation. These are defined in Section~\ref{spaces}.

Instead of working with the linear Schr\"odinger equation we need to
consider lower order perturbations of it of the form
\begin{equation} \label{e2} \left\{
    \begin{aligned} 
      & iv_{t}-\Delta v-a(t,x,D)v= f(x,t) \\& v(x,0) =g(x)
    \end{aligned}
  \right.
\end{equation}
Given any $M \geq 1$ we introduce a larger class $\mathcal{C}_M$ of
pseudodifferential operators so that for $a \in \mathcal{C}_M$ we can
solve \eqref{e2}.  To motivate the following definition, we note that
the imaginary part of the symbol $a$ can produce exponential growth in
\eqref{e2}. We want to be able to control the size of the growth in
the phase space along the Hamilton flow of the linear Schr\"odinger
equation, which leads to a (possibly large) bound on the integral of
$a$ along the flow. We also want the integral of $a$ along the flow to
give an accurate picture of the evolution. To insure this we impose a
smallness condition on the integral of derivatives of $a$ along the
flow. This motivates the following

\begin{d1} Let $M \geq 1$. The symbol $a(x,\xi)$ belongs to the class
  $\mathcal{C}_M$ if it satisfies the following conditions:

  \beq \label{aa1} M= \sup_{x,\xi} \int_{0}^1 |a(t,x+2t\xi,\xi)| dt <
  \infty \eeq

  \beq \label{aa2} \sup_{x,\xi} \int_{0}^1 |\partial_x^\alpha
  \partial_\xi^\beta a(t,x+2t\xi,\xi)| dt \leq C_{\alpha,\beta} \delta
  \qquad |\alpha|+|\beta| \geq 1 \eeq

\noindent
with $\delta e^{M} \ll 1$.
\end{d1}

Then our main linear result has the form:
 
\begin{t1} \label{le} Let $s \in \R$ and $a \in \mathcal{C}_M$. Then
  the solution $v$ of the equation \eqref{e2} satisfies the estimate:
  \begin{equation} \label{a2} \|\chi_{[0,1]}v\|_{\mathcal{D}X^{s}}
    \les e^{M}( \|f\|_{\mathcal{D}Y^{s}}+\|g\|_{\mathcal{D}H^{s}})
  \end{equation}
\end{t1}

To apply this theorem in our context we need to show that
\begin{p1} \label{low2} For $ u_0 \in \dh^s$ as above set
  \[
  a(x,\xi) = B_0((u_{0,\leq 1}^{\e},\bar u_{0,\leq 1}^{\e}),\xi)
  \]
  If $\e \leq e^{-CM}$ with $C$ sufficiently large then $a \in
  \mathcal{C}_{c_n M}$ with $c_n$ depending only on the dimension $n$.
\end{p1}

In order to iteratively solve the equation \eqref{e1} we want bounds
for the right hand side terms. The main one is a bilinear estimate,

\begin{t1} \label{ta3} The following bilinear estimate holds:
  \begin{equation}
    \label{a3} \|B((u,\bar u),(\nabla v,\nabla \bar v))\|_{\mathcal{D}Y^{s}} \lesssim 
    \|u\|_{\mathcal{D}X^{s}}\|v\|_{\mathcal{D}X^{s}} \qquad s > \frac{n}2
  \end{equation}
\end{t1}
For the linear term $A_1$ we use the estimates
\begin{p1} \label{lowlin} Let $A_1$ be defined as above. Then
  \begin{equation}
    \|A_1 v\|_{Y^s} \lesssim \| u_{0, \leq 1}^{\e}\|_{L^\infty} \|
    v\|_{X^s},
    \qquad
    \|A_1 v\|_{\dy^s} \lesssim \| u_{0, \leq 1}^{\e}\|_{L^\infty} \| v\|_{\dx^s}.
 \label{lowe} \end{equation}

\end{p1}

and one for the time independent term:

\begin{p1} \label{low3} If $s > 1$ then
  \begin{equation}
    \label{a6} \|B(u_{0, \leq 1}^{\e}, \nabla{u_{0, \leq 1}^{\e}})+ \DD
    u_{0, \leq 1}^{\e}\|_{\mathcal{D}H^{s}} \lesssim \min{(\e, \e^{s-\frac{n}2})} M^2
  \end{equation}
\end{p1}

Using the results above we are able to conclude the proof of
Theorem~\ref{maint}. We choose $\e = e^{-C M}$ with $C$ large enough
(depending on $s$). Then we solve the rescaled problem \eqref{e1} on a
unit time interval using the contraction principle in the space $X^s$.

We define the operator $\mathcal{T}_{1}$ by $w=\mathcal{T}_{1}f$ to be
the solution of the inhomogeneous Schr\"odinger equation with zero
initial data:
\begin{equation} \label{e3}
  \begin{cases}
    \begin{aligned}
      &iw_{t}-\Delta w -A(x,D) w = f \\
      &w(x,0) =0
    \end{aligned}
  \end{cases}
\end{equation}
We also denote by $\mathcal{T}_{2}g$ the solution to the homogeneous
equation:
\begin{equation} \label{e4}
  \begin{cases}
    \begin{aligned}
      &iw_{t}-\Delta w -A(x,D) w = 0 \\
      &w(x,0) =g
    \end{aligned}
  \end{cases}
\end{equation}
With these notations the equation \eqref{e1} can be rewritten in the
form
\[
v = \mathcal{T} v, \qquad \mathcal{T} v = \mathcal{T}_{2} u_{0, >
  1}^\e + \mathcal{T}_{1}(B(v,\nabla v) + A_1 v + B(w_{0}^{\e},
\nabla{w_{0}}^{\e})+ \DD w_{0}^{\e})
\]
 
We define the set
\[
K=\{w \in \mathcal{D}X^{s}: \|w\|_{\mathcal{D}X^{s}} \leq \e^\sigma
\}, \qquad 0 < \sigma < s - \frac{n}2
\]
and prove that $\mathcal{T}: K \rightarrow K$ and that $\mathcal{T}$
is a contraction on $K$. This give us the existence of a fixed point
for $\mathcal{T}$ which is the solution of our problem in the interval
$[0,1]$.

To prove the invariance of $K$ under the action of $\mathcal{T}$ we
use the results in Theorem \ref{le}:
\[
\begin{split}
  \|\mathcal{T} v\|_{\mathcal{D}X^{s}} \lesssim & \ e^{c_n M}(
  \|u_{0,> 1}^{\e}\|_{\mathcal{D}H^{s}} + \| B(v,\nabla
  v)\|_{\mathcal{D}Y^{s}} + \|A_1 v\|_{\mathcal{D}Y^{s}} \\ + &\ \|B(u_{0,\leq 1}^{\e}
  \nabla{u_{0,\leq 1}^{\e}})+ \DD u_{0,\leq
    1}^{\e}\|_{\mathcal{D}H^{s}})
\end{split}
\]
Using the estimates (\ref{uhigh}), (\ref{a3}), (\ref{lowe}) and (\ref{a6}) this
yields
\[
\|\mathcal{T}v\|_{\mathcal{D}X^{s}} \lesssim M^2 e^{c_n M}
(\e^{s-\frac{n}2} + \| v\|_{\mathcal{D}X^{s}}^{2} )
\]
If $v \in K$ then we use the smallness of $\e$ to obtain
\[
\|\mathcal{T}v\|_{\mathcal{D}X^{s}} \lesssim \e^{\sigma} M^2 e^{c_n M}
(\e^{s-\frac{n}2-\sigma} + \e^\sigma) \leq \e^{\sigma}
\]
which shows that $\mathcal{T}v \in K$.

To prove that $\mathcal{T}$ is a contraction we write
\[
\begin{split}
  \mathcal{T}v_{1}-\mathcal{T}v_{2} =\mathcal{T}_{1}(
  B(v_{1},\nabla(v_{1}-v_{2})))+\mathcal{T}_{1}( B(v_1-v_{2},\nabla
  v_{2}))
\end{split}
\]
and estimate in a similar manner
\[
\begin{split}
  \|\mathcal{T}v_{1}-\mathcal{T}v_{2}\|_{\mathcal{D}X^{s}} \les& (\|
  v_{1}\|_{\mathcal{D}X^{s}}+
  \|v_{2}\|_{\mathcal{D}X^{s}})\|v_{1}-v_{2}\|_{\mathcal{D}X^{s}} \\
  \les & 2 \e^\sigma \|\chi_{[0,1]}(v_{1}-v_{2})\|_{\mathcal{D}X^{s}}
  < \frac12 \|\chi_{[0,1]}(v_{1}-v_{2})\|_{\mathcal{D}X^{s}}
\end{split}
\]

\section{The function spaces} \label{spaces}

Let $(x_0,\xi_0)\in \R^{2n}$. To describe functions which are
localized in the phase space on the unit scale near $(x_0,\xi_0)$
we use the norm:
\[
H^{N,N}_{x^0,\xi^0}:=\{f: \langle D-\xi^0\rangle ^{N} f \in L^2,\ \langle x-x^0
\rangle^{N} f \in L^2 \}
\]


We work with the lattice $\Z^n$ both in the physical and Fourier space.
 We consider a partition of unity in the physical space,
\[
\sum_{x_0 \in \Z^n} \phi_{x_0}=1  \qquad \phi_{x_0}(x) = \phi(x-x_0)
\]
where $\phi$ is a smooth bump function with compact support.
We  use a similar partition of unity on the
Fourier side:
\[
\sum_{\xi_0 \in \Z^n} \varphi_{\xi_0} =1, \qquad \varphi_{\xi_0}(\xi) = \varphi(\xi-\xi_0)
\]

Let $H$ be a Hilbert space.  Let $V^{2}H$ be the space of continuous
$H$ valued functions on $\R$ with bounded 2-variation:
\[
\|u\|^{2}_{V^{2}H} =\sup_{(t_{i}) \in T} \sum_{i}
\|u(t_{i+1})-u(t_{i})\|_H^{2}
\]
where $T$ is the set of finite increasing sequences in $\R$. We abuse
notations and consider $V^{2}$ to be the closure of the space of
smooth functions in this norm, namely we use $V^{2}$ instead of $V^{2}
\cap C$. The $V^2$ spaces are close to the homogeneous Sobolev space
$\dot H^\frac12$ in the sense that
\begin{equation}
\dot B^{\frac12}_{2,1} \subset V^2 \subset \dot B^{\frac12}_{2,\infty}
\label{emb}\end{equation}

Let $U^{2}H$ be the atomic space defined by the atoms:
\[
u=\sum_{i} h_{i} \chi_{[t_{i},t_{i+1})}, \ \sum_{i} \|h_{i}\|_H^{2}=1
\]
\noindent
for some $(t_{i}) \in T$. We have the inclusion $U^{2}H \subset V^2H$
but in effect these spaces are very close.  There is also a duality
relation between $V^{2}H$ and $U^{2}H$, namely
\begin{equation}
(DU^{2}H)^{*}=V^{2}H
\label{dual1}\end{equation}

We can associate similar spaces associated to the Schrodinger flow
by pulling back functions to time $0$ along the flow, e.g.
\[
\| u \|_{V^2_\Delta L^2} = \| e^{it\DD} u\|_{V^2 L^2} 
\]
This turns out to be a good replacement for the $X^{0,\frac12}$ space
associated to the Schr\"odinger equations. Such spaces originate in 
unpublished work of the second author on the wave-map equation,
and have been succesfully used in various contexts so far, see
\cite{kt} and \cite{ss}. 

In the present paper we consider a wave packet type refinement of this
structure.  We begin with spatial localization, and introduce the
space $X$ of functions in $[0,1] \times \R^n$ with norm
\[
\| u\|_{X}^2 = \sum_{x_0 \in Z^2} 
\| \phi_{x_0} e^{it\DD} u\|_{V^2_t L^2_x}^2 
\]
It is easy to see that this is a stronger norm than the $V^2_\Delta
L^2$ norm,
\begin{equation}
\| u \|_{V^2_\Delta L^2} \lesssim \| u\|_{X}^2
\label{l2d} \end{equation}

The $X$ norm is usually applied to functions which are also frequency
localized on the unit scale. This is consistent with the spatial
localization. Precisely, we have the straightforward bound
\begin{equation}\label{xloc}
\| \varphi_\xi(D) u\|_{X} \lesssim \| u\|_{X}
\end{equation}

 For a function $u :[0,1] \times  R^{n} 
\rightarrow \C$ we decompose
\[
u=\sum_{\xi_0 \in \Z^n} u_{\xi_0}, \qquad  u_{\xi_0}={\varphi}_{\xi_0}(D)u
\]
and set
\[
\| u\|_{X^s}^2 = \sum_{\xi_0 \in \Z^n} \langle \xi_0 \rangle^{2s} 
\|u_{\xi_0} \|_{X}^2
\]
An immediate consequence of \eqref{l2d} is that
\begin{equation}
\| u \|_{V^2_\Delta L^2} \lesssim \| u\|_{X^0}^2
\label{x0x2}\end{equation}

We also denote by $X_\lambda$ the subspace of functions in $X^0$ which 
are localized at frequency $\lambda$. It is easy to see that $X^s$ has an
$l^2$  dyadic structure,
\[
\| u\|_{X^s}^2 \approx \sum_{\lambda \geq 1} \lambda^{2s} \| u_\lambda
\|_{X_\lambda}^2 
\]

This is the most compact definition of $X^s$. Some equivalent
formulations of the norm on $X$ turn out to be more helpful in
some estimates. We dedicate the next few paragraphs to such
formulations. The first such formulation simply adds regularity and
decay:

\begin{p1} \label{char22}
Let $N \in \N$. Then 
\beq \label{char1}
\| u_{\xi_0}\|_{X}^2 \approx \sum_{x_0 \in Z^2} 
\|  \phi_{x_0} e^{it\DD} u_{\xi_0}\|_{V^2 H^{N,N}_{x_0,\xi_0}}^2 
\eeq
\end{p1}

\begin{proof}
 Denote $v_{\xi_0} =  e^{it\DD} u_{\xi_0}$. Then it suffices to show
 that
\[
\sum_{x_0 \in Z^2} \| \phi_{x_0} v_{\xi_0}\|_{V^2 H^{N,N}_{x_0,\xi_0}}^2 
\lesssim \sum_{x_0 \in Z^2} \|  \phi_{x_0} v_{\xi_0}\|_{V^2 L^2}^2 
\]
This in turn follows by summation from 
\[
 \| \phi_{x_0} v_{\xi_0}\|_{V^2 H^{N,N}_{x_0,\xi_0}} \lesssim
\sum_{y_0 \in \Z^n} \langle x_0 -y_0 \rangle^{-N}
 \|  \phi_{y_0} v_{\xi_0}\|_{V^2 L^2}
\]
We can translate both in space and in frequency and reduce the problem
to the case when $x_0=0$ and $\xi_0=0$. Then we need to show that
\[
\| x^N \phi_0 v_0\|_{V^2 L^2} + \| D^N (\phi_0 v_0)\|_{V^2 L^2}
\lesssim \sum_{y_0 \in \Z^n} \langle y_0 \rangle^{-N} \| \phi_{y_0}
v_{0}\|_{V^2 L^2}
\]
The derivatives which fall on $v_0$ can be truncated at frequencies 
larger than $1$. Then a more general formulation of the above bound is
\begin{equation}
\| \psi_0 \chi_0(D) v_0\|_{V^2 L^2} \lesssim \sum_{y_0 \in \Z^n} 
\langle y_0 \rangle^{-N}  \|  \phi_{y_0} v_{0}\|_{V^2 L^2}
\label{offdecay}\end{equation}
where both $\psi_0$ and $\chi_0$ are bump functions concentrated at
$0$. We write
\[
 \psi_0 \chi_0(D)  v_0 = \sum_{y_0 \in \Z^n}   \psi_0 \chi_0(D) \phi_{y_0} v_{0}
\]
It remains to show that
\[
\|  \psi_0 \chi_0(D) \phi_{y_0}\|_{L^2 \to L^2} \lesssim \langle y_0
\rangle^{-N}
\]
which is straightforward since $\chi_0(D)$ has a bounded and rapidly
decreasing kernel.
\end{proof}

The next equivalent definition of our function spaces relates the
Schr\"odinder evolution to the associated Hamilton flow,
\[
(x_0,\xi_0) \to (x_0-2t\xi_0, \xi_0)
\]
We begin by the linearizing  the symbol of $-\Delta$ near
$\xi_0$,
\[
\xi^2 = L_{\xi_0}(\xi) + O((\xi-\xi_0)^2), \qquad 
 L_{\xi_0}(\xi) = \xi_0^2 + 2 \xi \xi_0 
\]
The evolution generated by $L$ is simply the transport along the 
Hamilton flow,
\[
e^{itL_{\xi_0}} u(x) = e^{it \xi_0^2} u(x-2t\xi_0)
\]
Our next characterization of the frequency localized $X$ norm asserts that we 
can replace $-\Delta$ by $L_{\xi_0}$.

\begin{p1}
We have
\begin{equation}
\| u_{\xi_0}\|_{X}^2 \approx  \sum_{x_0 \in Z^2} 
\|  \phi_{x_0} e^{-itL_{\xi_0}} u_{\xi_0}\|_{V^2 L^2}^2 
\label{eqwp}\end{equation}
and
\begin{equation}
\| u_{\xi_0}\|_{X}^2 \approx \sum_{x_0 \in Z^2} 
\| \phi_{x_0} e^{-it L_{\xi_0}} u_{\xi_0}\|_{V^2 H^{N,N}_{x_0,\xi_0}}^2 
\label{eqwpN}\end{equation}
\end{p1}
\begin{proof}

Denoting as before $v_{\xi_0} = e^{it\DD} u_{\xi_0}$ we can write
\[
e^{-itL_{\xi_0}} u_{\xi_0} = e^{it(D^2-L_{\xi_0})} v_{\xi_0} =
\chi_{\xi_0}(t,D) v_{\xi_0} 
\]
where the symbol $\chi(t,\xi)$ is a unit bump function in $\xi$ around
$\xi_0$ and smooth in $t$.  There is also a similar formula with
$L_{\xi_0}$ and $D^2$ interchanged. Hence for \eqref{eqwp} it suffices
to show that
\[
 \sum_{x_0 \in Z^2}  \|  \phi_{x_0}  \chi_{\xi_0}(t,D) v_{\xi_0}\|_{V^2 L^2}^2 
\lesssim  \sum_{x_0 \in Z^2}  \|  \phi_{x_0} v_{\xi_0}\|_{V^2 L^2}^2 
\]
Without any restriction in generality we can take $\xi_0=0$. Using a
Fourier series in $t$ we can also replace $\chi_{\xi_0}(t,\xi)$ by an
expression of the form $a(t) \chi_{\xi_0}(\xi)$ with $a$ smooth. It
remains to show that 
\[
 \sum_{x_0 \in Z^2}  \|  \phi_{x_0}  \chi_{\xi_0}(D) v_{\xi_0}\|_{V^2 L^2}^2 
\lesssim  \sum_{x_0 \in Z^2}  \|  \phi_{x_0} v_{\xi_0}\|_{V^2 L^2}^2 
\]
But this follows as in Proposition~\ref{char22} from \eqref{offdecay}.

Finally, the proof of \eqref{eqwpN} uses the same argument as in
Proposition~\ref{char1}.
\end{proof}

For each $\xi_0$ we define the set of tubes $T_{\xi_0}$ of the form
\[
 Q = \{  (t,x): |x-(x_Q-2t\xi_0)| \leq 1 \}  \qquad x_Q \in \Z^n
\]
They are the image of unit cubes centered at $x_Q$ at time $0$ 
along the $L_{\xi_0}$ flow.  These tubes generate the decomposition:  
\[
[0,1] \times \R^n = \bigcup_{Q \in T_{\xi_0}} Q
\]

Now we can state our last equivalent formulation of the $X$ norm for
functions localized at frequency $\xi_0$ in terms of a wave packet
decomposition associated to tubes $Q \in T_{\xi_0}$:

\begin{p1}
Let $u_{\xi_0} \in X$. Then it can be represented as the sum a rapidly
convergent series
\begin{equation}
u_{\xi_0}(t,x) =   \sum_j  e^{ix\xi_0} e^{it\xi_0^2} \sum_{Q \in T_{\xi_0}}  
a_Q^j(t) \chi^j(x-x_Q - 2t\xi_0)
\end{equation} 
with $\chi^j$ uniformly bounded in $H^{N,N}$, supported in $B(0,2)$ and
\[
\sum_{Q \in T_{\xi_0}}  
\|a_Q^j\|_{V^2}^2 \lesssim j^{-N} \|u_{\xi_0}\|_{X}
\]
\label{rep}\end{p1}

\begin{proof}
The pull back to time $0$ of the above representation using the
$L_{\xi_0}$ flow is
\[
e^{-itL_{\xi_0}} u_{\xi_0}(t,x) = \sum_j  e^{ix\xi_0}  \sum_{x_0 \in \Z^n}   
a_Q^j(t) \chi^j(x-x_0)
\]
We denote $v_{0} = e^{-ix\xi_0} e^{-itL_{\xi_0}} u_{\xi_0}$ and 
use  \eqref{eqwpN} for the $X_{\xi_0}$ norm.
Then it suffices  to show that for fixed $x_0$ we can represent
\[
\phi_{x_0} v_0 =  \sum_j  a^j(t) \chi^j(x-x_0)
\]
where $ \chi^j$ are uniformly bounded in $H^{N,N}$ and
\[
\|a^j\|_{V^2} \lesssim j^{-N} \|\phi_{x_0} v_0 \|_{V^2 H^{CN,CN}_{x_0}}
\]
Without any restriction in generality we take $x_0=0$. We take
$\chi^j$ to be the Hermite functions with the $H^{N,N}$ normalization
(see e.g. \cite{Taylor}).  Then the $V^2$ functions $a^j$ are the
Fourier coefficients of $\phi_{x_0}v_0$. They decay rapidly due to the
additional regularity of $\phi_{x_0} v_0$.

Finally, to insure that the $\chi_j$'s have compact support we can truncate the
Fourier series outside the support of $\phi_{x_0}$. 

\end{proof}

We follow a similar path to define the $Y^{s}$ structure:
\[
\| u\|_{Y^s}^2 = \sum_{\xi_0 \in \Z^n} \langle \xi_0 \rangle^{2s} 
\|u_{\xi_0} \|_{Y}^2
\]
where $Y$ is defined by:
\[
\|f \|_{Y}^2 = \sum_{x_0 \in Z^2} \| \phi_{x_0} e^{-itD^2} f\|_{DU^2_t L^2_x}^2. 
\]

All the equivalent definitions for the $X$ have a counterpart
for $Y$ by simply replacing the $V^2_t$ structure by $DU^2_t$ one.

There are two key relations between the $X^s$ and the $Y^s$ spaces.
The first one is concerned with solvability for the linear
Schr\"odinger equation:

\begin{p1}
The solution $u$ to the linear Schr\"odinger equation
\[
iu_t -\Delta u = f, \qquad u(0) = u_0
\]
satisfies
\begin{equation}
\| u\|_{X^s} \lesssim \|u_0\|_{H^s} + \|f\|_{Y^s}
\end{equation}
\end{p1}

This is stated here only for the sake of completeness, as in the next
section we prove the stronger estimate in Theorem~\ref{le}. 

The second relation is a duality relation:
\begin{p1}
We have the duality relation
\[
(Y^s)^* = X^{-s}.
\]
\label{pdual}\end{p1}

\begin{proof}
a) We first verify that $X^{-s} \subset (Y^s)^*$. For this we need the bound
\[
\left | \int_{0}^1 \int u \bar f dx dt\right | \lesssim \|u\|_{X^{-s}} \|f\|_{Y^s}
\]
We decompose 
\[
\int _{0}^1 \int u \bar f dx dt = \sum_{\xi\in \Z^n} \int _{0}^1 \int
u_\xi \bar f_\xi dx dt
\]
Due to the definition of the $X^s$ and $Y^s$ norms, it suffices to show
that 
\[
\left |\int _{0}^1 \int u_\xi \bar f_\xi dx dt\right|   
\lesssim \|u_\xi\|_{X} \|f_\xi\|_{Y}.
\]
For this we write
\[
\begin{split}
\int _{0}^1 \int u_\xi \bar f_\xi dx dt =& 
\int _{0}^1 \int e^{it\DD} u_\xi \ \overline{e^{it\DD} f_\xi} dx dt \\
=& \sum_{x_0 \in Z^2} \int _{0}^1 \int \phi_{x_0} 
e^{it\DD} u_\xi \ \overline{e^{it\DD} f_\xi} dx dt
\end{split}
\]
and use the duality relation \eqref{dual1} together with the definition
of the $X$ and $Y$ norms.

b) We now show that $ (Y^s)^* \subset X^{-s}$.
Let $T$ be a bounded linear functional. Then we have
\[
\begin{split}
|T f| \lesssim \|f\|_{Y^s} &\lesssim \left( \sum_{\xi \in \Z^n} 
\langle \xi \rangle^{2s} \|f_\xi\|_{Y}^2 \right)^{\frac12} 
\\ & \lesssim \left( \sum_{x_0 \in Z^2} \sum_{\xi \in \Z^n} 
\langle \xi \rangle^{2s} \|\phi_{x_0} e^{it\DD} f_\xi\|_{DU^2L^2}^2 
\right)^{\frac12} 
\end{split}
\]
By the Hahn-Banach theorem we can extend $T$ to a bounded 
linear functional on the space $l^2_{ \langle \xi \rangle^{s}} DU^2L^2$.
By \eqref{dual1} this implies that we can represent $T$ in the form
\[
Tf = \sum_{x_0 \in Z^2} \sum_{\xi \in \Z^n}  \int _{0}^1 
\int v_{x_0,\xi}  \overline{\phi_{x_0} e^{it\DD} f_\xi}  dx dt 
\]
where 
\[
\sum_{x_0 \in Z^2} \sum_{\xi \in \Z^n}\langle \xi \rangle^{-2s}
\|v_{x_0,\xi} \|_{V^2 L^2}^2 \lesssim \|T\|_{ (Y^s)^*}^2
\]
Due to the above representation we can identify $T$ with the function 
\[
u_T =  \sum_{x_0 \in Z^2} \sum_{\xi \in \Z^n}   
 \phi_\xi(D) e^{-it\DD} \phi_{x_0} v_{x_0,\xi}
\]
It remains to show that 
\[
\|u_T\|_{X^{-s}}^2 \lesssim \sum_{x_0 \in Z^2} \sum_{\xi \in
  \Z^n}\langle \xi \rangle^{-2s} \|v_{x_0,\xi} \|_{V^2 L^2}^2
\]
This reduces to the fixed $\xi$ bound
\[
\|   \sum_{x_0 \in Z^2} \varphi_\xi(D) e^{-it\DD} \phi_{x_0}
v_{x_0,\xi}\|_{X}^2 \lesssim  \sum_{x_0 \in Z^2}  \|v_{x_0,\xi} \|_{V^2 L^2}^2
\]
with a modified $\varphi_\xi$. Using the definition of the $X$ norm
we rewrite this as 
\[
\sum_{y_0 \in Z^2} \| \sum_{x_0 \in Z^2} \phi_{y_0} \varphi_\xi(D)
 \phi_{x_0} v_{x_0,\xi}\|_{V^2 L^2}^2
\lesssim  \sum_{x_0 \in Z^2}  \|v_{x_0,\xi} \|_{V^2 L^2}^2.
\]
But this follows by Cauchy-Schwartz from the rapid decay
\[
\|  \phi_{y_0} \varphi_\xi(D)
 \phi_{x_0}\|_{L^2 \to L^2} \lesssim \langle x_0-y_0\rangle^{-N}.
\]
\end{proof}

We still need to add the decay structure to the $X^s$, respectively
the $Y^s$ spaces. Given  $\lambda \geq 1$, we roughly want to ask
for $l^1$ summability of frequency $\lambda$ norms along collinear cubes
of size $\lambda$.  We define the 
$\dx_\lambda$ norm by
\begin{equation}
\| u\|_{\dx_\lambda} = \sup_{x_0 \in \R^n} \sup_{\omega \in \S^{n-1}}
\sum_{k = -\infty}^\infty \| \chi(\lambda^{-1}(x-x_0) - k \omega) u\|_{X_{\lambda}} 
\end{equation}
where $\chi$ is a compactly supported bump function.  We note that
this norm can only be meaningfully used for functions at frequency
$\lambda$. Summing up with respect to $\lambda$ we also set
\[
\| u\|_{\dx^s}^2 = \sum_{\lambda \geq 1} \|  u_\lambda\|_{\dx_\lambda}^2.
\]
In a completely similar way we can define $\mathcal{D}Y_{\lambda}$ and $\mathcal{D}Y^{s}$.

A useful tool in our analysis is a family of embeddings which
correspond to the Strichartz estimates for the Schr\"odinger equation.

\begin{p1}
Let $p$ and $q$ be indices which satisfy
\[ 
\frac{2}p+ \frac{n}q = \frac{n}2, 2 < p \leq \infty, \ 2 \leq  q \leq
\infty
\]
Then we have the embedding 
\[
X^0 \subset L^p_t L^q_x
\]
\label{SE}\end{p1} 

By \eqref{x0x2} these embeddings are a direct consequence of the
Strichartz estimates for the Schr\"odinger equation, see \cite{kt},
the proof of Proposition 6.2.  By Sobolev embeddings we also obtain
bounds with larger $p,q$ for frequency localized solutions.

\begin{c1}
Let $p$ and $q$ be indices which satisfy
\[ 
\frac{2}p+ \frac{n}q \leq  \frac{n}2,\qquad  2 < p \leq \infty,\qquad  2 \leq  q \leq
\infty
\]
Then we have the embedding 
\[
S_{< \mu}  X^0 \subset \mu^{ \frac{n}2 - \frac{2}p- \frac{n}q}  L^p_t L^q_x
\]
where $S_{<\mu}$ can be replaced with a multiplier localizing in
frequency to an arbitrary cube of size $\mu$.
\end{c1}

In this article we  use only the case $p=q$, more precisely
\begin{equation}
 X^0  \subset L^\frac{2(n+2)}{n}, \qquad 
S_{<\mu} X^0  \subset \mu^{\frac{n}2-1} L^{n+2}
\label{set}\end{equation}

Finally, we introduce modulation localization operators, which we can
define in two equivalent ways. First is as multipliers,
\[
\widehat{M_{<\sigma} u} = s_{<\sigma}(\tau-\xi^2) \hat u.
\]
The second is obtained by conjugation with respect to the 
Schr\"odinger flow,
\[
 e^{itD^2} (M_{<\sigma} u)(t) = s_{<\sigma}(D_t) (e^{itD^2} u(t))  
\]

The operators $s_{<\sigma}(D_t)$ are bounded on $V^2$. On the other
hand  for the remainder, by \eqref{emb},  we have a good $L^2$ bound:
\[
\| s_{>\sigma}(D_t) a\|_{L^2} \lesssim \sigma^{-\frac12} \|a\|_{V^2}.
\]
By the second definition of the modulation localization operators 
above we obtain

\begin{p1}
a) The operators $M_{<\sigma}$ are bounded on $X^0$.

b) The following estimate holds:
\[
\| M_{>\sigma} u\|_{L^2} \lesssim \sigma^{-\frac12} \|u\|_{X^0}
\]
\label{mod}\end{p1}

We note that in  effect, by using both inclusions in \eqref{emb}, one
can relate our $X^s$ spaces to the traditional $X^{s,b}$ spaces,
namely
\begin{equation}
\dot{X}^{s,\frac12,1} \subset X^s \subset \dot{X}^{s,\frac12,\infty}
\end{equation}

\section{Linear estimates}
\label{linearestimates}

This section is devoted to the study of the linear equation
\eqref{e2} with $ a \in \mathcal{C}_M$.  Precisely,  we aim to
prove Theorem~\ref{le}. We denote
\[
L_a = i \partial_t - \Delta - a(t,x,D)
\]

We begin our analysis with a heuristic computation.  Suppose we have a
solution $u^{x_0,\xi_0}$ to
\[
(i \partial_t - \Delta - a(t,x,D) ) u = 0
\]
which is localized on the unit scale near the bicharacteristic $t \to
(x_0+2t\xi_0,\xi_0)$. Then we can freeze the symbol of $a(t,x,D)$
along the ray and write
\[
a(t,x,D) u = a(t,x_0+2\xi_0,\xi_0) u + error
\]
Thus $u^{x_0,\xi_0}$ approximately solves
\[
(i \partial_t - \Delta - a(t,x_0+2\xi_0, \xi_0)) u_{x_0,\xi_0} \approx
0
\]
This implies that we can represent $u^{x_0,\xi_0}$ in the form
\[
u^{x_0,\xi_0}(t) \approx e^{\int_{0}^t a(s,x_0+2s\xi_0, \xi_0) ds}
v^{x_0,\xi_0}
\]
where $v^{x_0,\xi_0}$ solves the equation
\[
(i \partial_t - \Delta ) v^{x_0,\xi_0} = 0
\]
Hence along each wave packet we can use the above exponential to
approximately conjugate the variable coefficient equation to the flat
flow.

Using this idea we produce wave packet approximate solutions for the
equation \eqref{e2}. By orthogonality these combine into general
approximate solutions to \eqref{e2}. The exact solutions are obtained
via a Picard iteration.

We first consider the regularity of the exponential weight. By
\eqref{aa1} we have
\[
\int_{0}^1 | a(t,x_0+2t\xi_0, \xi_0) |dt \leq M
\]
which implies a $W^{1,1}$ bound for the exponential,
\begin{equation}
  \left \| \frac{d}{dt} e^{\int_{0}^t a(t,x_0+2t\xi_0, \xi_0) dt} \right \|_{L^1}
  \leq e^M
  \label{wiiint}\end{equation}
Since $W^{1,1} \subset V^2$ this yields a similar $V^2$ bound,
\begin{equation}
  \left \| e^{\int_{0}^t a(t,x_0+2t\xi_0, \xi_0) dt} \right \|_{V^2}
  \leq e^M
  \label{v2int}\end{equation}

We continue our analysis with the localized equation
\begin{equation}
  L_a u =  f_{x_0,\xi_0}, \qquad
  u(0) =  g_{x_0,\xi_0}
  \label{x0xi0}\end{equation}
According to the above heuristics, an approximate solution for this
should be given by Duhamel's formula,
\[
\begin{split}
  u^{x_0,\xi_0} =& \ \ e^{\int_{0}^t a(s,x_0+2s\xi_0, \xi_0 ds)}
  e^{-it\Delta} g_{x_0,\xi_0} \\ & + \int_0^t e^{\int_{s}^t a(\tau,
    x_0+2\tau \xi_0, \xi_0) d \tau} e^{-i(t-s)\Delta} f_{x_0,\xi_0}(s)
  ds
\end{split}
\]
We prove that this is indeed the case:

\begin{p1} 
  The function $ u^{x_0,\xi_0} $ is an approximate solution for
  \eqref{x0xi0} in the sense that
  \begin{equation} \label{sol} \|e^{it\Delta} u^{x_0,\xi_0}\|_{V^2
      H^{N,N}_{x_0,\xi_0}} \lesssim e^{M}
    (\|g_{x_0,\xi_0}\|_{H^{N,N}_{x_0}}+ \| e^{it\Delta}
    f_{x_0,\xi_0}\|_{DU^2 H^{N,N}_{x_0}})
  \end{equation}
  and
  \begin{equation} \label{eerr}
    \begin{split} 
     & \| e^{it\Delta}\left( L_a
        u^{x_0,\xi_0}-f_{x_0,\xi_0} \right) \|_{L^1
        H^{N,N}_{x_0,\xi_0}} \lesssim \\
      & \delta e^{M} (\|g_{x_0,\xi_0}\|_{H^{N+2n+1,N+2n+1}_{x_0,\xi_0}}+ \|
      e^{it\Delta} f_{x_0,\xi_0}\|_{DV^2 H^{N+2n+1,N+2n+1}_{x_0,\xi_0}})
    \end{split}
  \end{equation}
\end{p1}

\begin{proof}
  For the first bound we shorten the notation
  \[
  a_0(s) := a(s,x_0+2s\xi_0, \xi_0)
  \]
  and compute
  \[
  e^{it \Delta} u^{x_0,\xi_0} = e^{\int_{0}^t a_0(s) ds} g_{x_0,\xi_0}
  + \int_0^t e^{\int_{s}^t a_0(\tau) d \tau} e^{is\Delta}
  f_{x_0,\xi_0}(s) ds
  \]
  The first factor is estimated directly by \eqref{v2int}. Setting
  \[
  F(t) = \int_0^t e^{is\Delta} f_{x_0,\xi_0}(s) ds
  \]
  and integrating by parts we write the second term in the form
  \[
  F(t) - \int_0^t \frac{d}{ds} e^{\int_{s}^t a(_0\tau) d \tau} F(s) ds
  \]
  We have
  \[
  \|F\|_{V^2H^{N,N}_{x_0,\xi_0}} \lesssim \| e^{is\Delta}
  f_{x_0,\xi_0}(s)\|_{DU^2 H^{N,N}_{x_0,\xi_0}}
  \]
  while the $V^2$ norm of the second part is controlled by its
  $W^{1,1}$ norm, namely
  \[
  \begin{split}
    \left\|\int_0^t \frac{d}{ds} e^{\int_{s}^t a(_0\tau) d \tau} F(s)
      ds \right\|_{V^2H^{N,N}_{x_0,\xi_0}} \lesssim \left \|
      \frac{d}{dt} \int_0^t \frac{d}{ds} e^{\int_{s}^t a_0(\tau) d
        \tau} F(s) ds\right\|_{L^1H^{N,N}_{x_0,\xi_0}} \hspace{-4in} &
    \\ \lesssim & \int_0^1 | a_0(t) |\left(1+\int_{0}^t \left|
        \frac{d}{ds} e^{\int_{s}^t a_0(\tau) d \tau}\right| ds\right)
    dt \|F\|_{L^\infty H^{N,N}_{x_0,\xi_0}} \\ \lesssim & \int_0^1
    |a_0(t) |\left(1+ \int_0^t |a_0(s)| e^{\int_{s}^t |a_0(\tau)|
        d\tau} ds \right) dt \|F\|_{L^\infty H^{N,N}_{x_0,\xi_0}} \\
    \lesssim & \int_0^1 |a_0(t)| e^{\int_{0}^t |a_0(\tau)| d\tau} dt
    \|F\|_{L^\infty H^{N,N}_{x_0,\xi_0}} \\ \lesssim &\
    e^{M}\|F\|_{L^\infty H^{N,N}_{x_0,\xi_0}} \lesssim \
    e^{M}\|F\|_{V^2 H^{N,N}_{x_0,\xi_0}}
  \end{split}
  \]
  This concludes the proof of \eqref{sol}.

  It remains to prove \eqref{eerr}. A direct computation yields
  \[
  (i\partial_t- \Delta - a(t,x,D) ) u^{x_0,\xi_0} -f_{x_0,\xi_0}
  =b(t,x,D) u^{x_0,\xi_0}
  \]
  where
  \[
  b(t,x,\xi) = a(t, x_0+2t\xi_0, \xi_0) - a(t,x,\xi)
  \]
  By \eqref{sol} it suffices to show that
  \[
  \| e^{it\DD} b(t,x,D) u^{x_0,\xi_0}\|_{L^1 H^{N,N}_{x_0,\xi_0}}
  \lesssim \delta \|e^{it\DD} u^{x_0,\xi_0}\|_{L^\infty
    H^{N+2n+1,N+2n+1}_{x_0,\xi_0}}
  \]
  Since the flat Schr\"odinger flow has the mapping property
  \[
  \| e^{it\DD} f\|_{H^{N,N}_{x_0,\xi_0}} \approx \|
  f\|_{H^{N,N}_{x_0+2t\xi_0,\xi_0}}
  \]
  this is equivalent to
  \begin{equation} \label{hna}
    \|  b(t,x,D) 
    u^{x_0,\xi_0}\|_{L^1 H^{N,N}_{x_0+2t\xi_0,\xi_0}} \lesssim \delta
    \| u^{x_0,\xi_0}\|_{L^\infty H^{N+3,N+3}_{x_0+2t\xi_0,\xi_0}}
  \end{equation}

  We begin with a straightforward consequence of the $S_{00}$
  calculus:

\begin{l1}
  Let $k$ be a nonnegative integer and $c$ be a symbol which satisfies
  \[
  |\partial_x^\alpha \partial_\xi^\beta c(x,\xi)| \leq
  c_{\alpha,\beta} \langle( x-x_0,\xi-\xi_0)\rangle^k, \qquad
  |\alpha|+|\beta| \geq 0
  \]
  Then for all $N \in \N$ we have
  \[
  \| c(x,D) u \|_{H^{N,N}_{x_0,\xi_0}} \lesssim
  \|u\|_{H^{N+k,N+k}_{x_0,\xi_0}}
  \]
  \label{soo}\end{l1}

In order to use this lemma for the operator $b$ above we need the
following Sobolev embedding:

\begin{l1}
  Let $R > 1$ and $c \in {W}^{2n,1}(B_R(x_0,\xi_0))$. Then
  \[
  |c(x,\xi)| \lesssim \sum_{0 \leq |\alpha| \leq 2n} \|\partial^\alpha
  c\|_{L^1(B_R(x_0,\xi_0))}
  \]
  respectively
  \[
  |c(x,\xi)-c(x_0,\xi_0)| \lesssim \sum_{1 \leq |\alpha| \leq 2n}
  \|\partial^\alpha c\|_{L^1(B_R(x_0,\xi_0))}
  \]
\end{l1}

As a consequence of this we obtain

\begin{l1}
  Let $c \in W_{loc}^{2n,1}(\R^{2n})$. Then
  \[
  \frac{|c(x,\xi)|}{\langle(x-x_0,\xi-\xi_0)\rangle^{2n+1} } \lesssim
  \int_{\R^{2n}} \sum_{0 \leq |\alpha| \leq 2n} \frac{
    |\partial^\alpha c(x,\xi)|}{\langle(x-x_0,\xi-\xi_0)\rangle^{2n+1}} dx
  d\xi
  \]
  respectively
  \[
  \frac{|c(x,\xi)-c(x_0,\xi_0)|}{\langle(x-x_0,\xi-\xi_0)\rangle^{2n+1}}
  \lesssim \int_{\R^{2n}} \sum_{1 \leq |\alpha| \leq 2n} \frac{
    |\partial^\alpha c(x,\xi)|}{\langle(x-x_0,\xi-\xi_0)\rangle^{2n+1}} dx
  d\xi
  \]
\end{l1}

Applying these inequalities to the symbol $b$ above we obtain
pointwise bounds for $b$
\[
\frac{|b(t,x,\xi)|}{\langle(x-x_0-2t\xi_0,\xi-\xi_0)\rangle^{2n+1}}
\lesssim \!\!
\int_{\R^{2n}}\!\!  \! \frac{ \sum_{ |\alpha|=1}^{2n} |\partial^\alpha
  a(t,x,\xi)|}{\langle(x-x_0-2t\xi_0,\xi-\xi_0)\rangle^{2n+1}} dx d\xi
\]
and also for its derivatives,
\[
\frac{| \partial^k b(t,x,\xi)|}{\langle(x-x_0-2t\xi_0,\xi-\xi_0)\rangle^{2n+1}
} \lesssim \int_{\R^{2n}}  \frac{\sum_{ |\alpha|=k}^{2n+k}
  |\partial^\alpha
  a(t,x,\xi)|}{\langle(x-x_0-2t\xi_0,\xi-\xi_0)\rangle^{2n+1}} dx d\xi.
\]
Then by Lemma~\ref{soo} it follows that $b$ satisfies the bound
\[
\| b(t,x,D)\|_{H^{N+2n+1,N+2n+1}_{x_0+2t\xi_0,\xi_0} \to
  H^{N,N}_{x_0+2t\xi_0,\xi_0}} \lesssim \int_{\R^{2n}}  \frac{ \sum_{
  |\alpha|=1}^{N_0}|\partial^\alpha
  a(t,x,\xi)|}{\langle(x-x_0,\xi-\xi_0)\rangle^{2n+1}} dx d\xi.
\]
with $N_0$ sufficiently large. Integrating in $t$ and changing
coordinates in the integral this gives
\[
\begin{split}
\| b(t,x,D)\|_{L^\infty H^{N+2n+1,N+2n+1}_{x_0+2t\xi_0,\xi_0} \!\! \to L^1
  H^{N,N}_{x_0+2t\xi_0,\xi_0}} \!\! \lesssim 
\\ \!\! \int_{\R^{2n}} \!\! \int_{0}^1
\sum_{|\alpha|=1}^{N_0} \frac{ |\partial^\alpha
  a(t,x+2t\xi,\xi)|}{\langle(x-x_0,\xi-\xi_0)\rangle^{2n+1}} dt dx
d\xi.
\end{split}
\]
By \eqref{aa2} we can bound each time integral by $\delta$ and the
remaining weight is integrable in $x$ and $\xi$. Hence \eqref{hna}
follows.

\end{proof}

Next we produce approximate solutions for the frequency localized
data,
\begin{equation}
  L_a u = f_{\xi_0}, \qquad u(0) = g_{\xi_0}
  \label{uxi}\end{equation}
We denote the approximate solution by $u^{\xi_0}$; the notation 
$u_{\xi_0}$ continues to be reserved for a frequency localized part 
of a function $u$.

\begin{p1} \label{luxi} There is an approximate solution $u^{{\xi_0}}$
  to the equation \eqref{uxi}, localized at frequency $\xi_0$, with
  $u(0) = g_{\xi_0}$, which satisfies the bounds
  \begin{equation} \label{sol2} \| u^{\xi_0} \|_{X} \lesssim
    e^{M} ( \|g_{\xi_0}\|_{L^2}+ \| f_{\xi_0}\|_{Y} )
  \end{equation}
  respectively
  \begin{equation} \label{eerr2} \| S_{\xi} ( L_a  u^{\xi_0} - f_{\xi_0} ) \|_{Y} \lesssim
    \delta \langle \xi -\xi_0 \rangle^{-N} e^{M} (
    \|g_{\xi_0}\|_{L^2}+ \| f_{\xi_0}\|_{Y} )
  \end{equation}
\end{p1}

\begin{proof}
  We decompose
\[
g_{\xi_{0}}=\sum_{x_{0}} g_{x_{0},\xi_{0}} \ \ \
f_{\xi_{0}}=\sum_{x_{0}} f_{x_{0},\xi_{0}}
\]
and solve the problem (\ref{x0xi0}) for which we have the estimates
(\ref{sol}) and (\ref{eerr}). We define our approximate solution to be
the sum of the approximate solutions $u^{x_0,{\xi_0}}$:
\[
u^{\xi_0}=\sum_{x_0} u^{x_0,\xi_0}
\]

Using (\ref{sol}) we obtain:
\[
\begin{split}
\|u^{\xi_0}\|^2_{X}=&\sum_{y_0} \|\phi_{y_0}
e^{-it\Delta}u_{\xi_0}\|^2_{V^2L^{2}} 
\\ \lesssim&
 \sum_{y_0} \left( \sum_{x_0} \|\phi_{y_0}
   e^{-it\Delta}u^{x_0,\xi_0}\|_{V^2L^{2}} \right)^2 
\\ \lesssim&
 \sum_{y_0} \left( \sum_{x_0} \langle x_0-y_0\rangle^{-N}
\|  e^{-it\Delta} u^{x_0,\xi_0}\|_{V^2 H^{N,N}_{x_0,\xi_0}} \right)^2 
\\ \lesssim&
  \sum_{x_0} \| e^{-it\Delta}u^{x_0,\xi_0}\|_{V^2L^{N,N}_{x_0,\xi_0}}^2
\\\lesssim&
 e^{M} \sum_{x_0} \|g_{x_0,\xi_0}\|_{H^{N+2n+1,N+2n+1}_{x_0,\xi_0}}^2+ 
\| e^{-it\Delta} f_{x_0,\xi_0}\|^2_{DU^2
  H^{N+2n+1,N+2n+1}_{x_0,\xi_0}}
\\ \lesssim&
 e^{M}
(\|g_{\xi_0}\|_{L^{2}}^2+ \| f_{\xi_0}\|^2_{Y_{\xi_0}})
\end{split}
\]

We continue now with the estimates for the error. using \eqref{eerr}
instead of \eqref{sol}.  We have
\[
\begin{split}
  \|S_{\xi} (L_a u^{\xi_0} -
  f_{\xi_0})\|_{Y}^2 =& \sum_y \|\phi_y e^{it\DD} S_{\xi} (L_a u^{\xi_0} -
  f_{\xi_0})\|_{DU^2 L^2}^2
  \\\lesssim&
  \sum_y  \left( \sum_{x_0}  
\| \phi_y S_{\xi}e^{it\DD} (L_a u^{x_0,\xi_0} - f_{x_0,\xi_0})\|_{L^1L^2}\right)^2
  \\\lesssim&   \sum_y  \left( \sum_{x_0}\frac{ 
\|e^{it\DD}( L_a u^{x_0,\xi_0} -
f_{x_0,\xi_0})\|_{L^1H^{N,N}_{x_0,\xi_0}}}{\langle(y-x_0-2t\xi_0,\xi-\xi_0)\rangle^{N}
}\right)^2
\\\lesssim&   \langle\xi-\xi_0\rangle^{n-2N} \sum_{x_0} 
\|e^{it\DD}( L_a u^{x_0,\xi_0} -
f_{x_0,\xi_0})\|_{L^1H^{N,N}_{x_0,\xi_0}}^2
\end{split}
\]
Then \eqref{eerr2} follows from \eqref{eerr}.
\end{proof}

The next stage is to consider data which is localized at frequency
$\lambda$,
\begin{equation}
  L_a u = f_{\lambda}, \qquad u(0) = g_{\lambda}
 \label{uj}\end{equation}
Summing up the frequency localized solutions we obtain as above
a dyadic approximate solution:

\begin{p1} \label{lulambda} There is an approximate solution $u^{\lambda}$
  to the equation \eqref{uj}, localized at frequency $\lambda$, with
  $u(0) = g_\lambda$, which satisfies the bounds
  \begin{equation} \label{sol3} 
\| u^\lambda \|_{X_\lambda} \lesssim
    e^{M} ( \|g_{\lambda}\|_{L^2}+ \| f_{\lambda}\|_{Y_\lambda} )
  \end{equation}
  respectively
  \begin{equation} \label{eerr3} \| S_\mu ( L_a u^{\lambda} -
    f_{\lambda} ) \|_{Y_\mu} \lesssim
    \left(\min\{\frac{\mu}\lambda,\frac\lambda\mu\}\right)^N \delta
    e^{M} (\|g_{\lambda}\|_{L^2}+ \| f_{\lambda}\|_{Y_j} )
  \end{equation}
\end{p1}

The construction of the functions $u^\lambda$ involves only the
constant coefficient Schr\"odinger flow at frequency $\lambda$. This
has spatial speed of propagation $\lambda$, therefore it can spread by
at most $O(\lambda)$ in a unit time interval. Thus the above
construction can be trivially localized on the $\lambda$ spatial
scale, leading to the bounds
 \begin{equation} \label{sol4} \| u^\lambda \|_{\dx_\lambda} \lesssim
    e^{M} ( \|g_{\lambda}\|_{\dl^2_\lambda}+ \| f_{\lambda}\|_{\dy_\lambda} )
  \end{equation}
  respectively
  \begin{equation} \label{eerr4} \| S_{\mu} ( L_a  u^{\lambda} - f_{\lambda} ) \|_{\dy_\mu} \lesssim
    \left(\min\{\frac{\mu}\lambda,\frac\lambda\mu\}\right)^N
 \delta e^{M} (\|g_{\lambda}\|_{\dl^2_\lambda}+ \| f_{\lambda}\|_{\dy_\lambda} )
  \end{equation}
After an addition dyadic summation we obtain a global
parametrix:

\begin{p1} \label{lu} There is an approximate solution $u$
  to the equation \eqref{e2}, $with u(0)=g$, which satisfies the bounds
  \begin{equation} 
\label{sol8} \| u \|_{\dx^s} \lesssim
    e^{M} ( \|g\|_{\dh^s}+ \| f\|_{\dy^s} )
  \end{equation}
  respectively
  \begin{equation} \label{eerr8} \| (i\partial_t
    -\Delta-a(t,x,D)) u - f ) \|_{\dy^s} \lesssim
    \delta  e^{M} (
    \|g\|_{\dh^s}+ \| f\|_{\dy^s} )
  \end{equation}
\end{p1}
Of course the similar result without decay is also valid.

If $\delta \ll e^{-M}$ then the constant in \eqref{eerr8} is less than
$1$. Then one can iterate to obtain an exact solution $u$ to
\eqref{e2} which still satisfies \eqref{sol4}. Theorem~\ref{le} follows.

\section{Low frequency bounds}

Here we prove Propositions~\ref{low1},\ref{low2},\ref{lowlin},\ref{low3}.

\begin{proof}[Proof of Proposition \ref{low1}] For simplicity we
  denote $f=u_0$.  The pointwise bounds \eqref{i2} follow from Sobolev type
  estimates and scaling,
\[
\| f^\e\|_{L^\infty} = \e \|f\|_{L^\infty} \lesssim \e
\|f\|_{\mathcal{D} H^s}, \qquad s > \frac{n}2
\]

We now move to the $L^2$ bounds. The effect of scaling on the
$\dl^2_\lambda$ norms is easy to compute,
\begin{equation}
\| f^\epsilon\|_{\dl^2_\lambda} = \e^{1-\frac{n}2} \| f\|_{\dl^2_{\e\lambda}} 
\label{dscale}\end{equation}
Hence the first bound in \eqref{ulow} can be rewritten in the form
\begin{equation}
 \e^{1-\frac{n}2}\| f_{<\e^{-1}}\|_{\dl^2_{\e}} \lesssim \|f\|_{\mathcal{D} H^s}
\label{7a}\end{equation}

We have 
  \[
 f_{\leq \e^{-1}}=\sum_{1 \leq \lambda \leq
  \e^{-1}} f_{\lambda} 
\]
For each such $\lambda$ we have three relevant spatial scales,
\[
\e \leq \lambda^{-1} \leq \lambda
\]
The middle one arises due to the uncertainty principle; namely, this
is the scale on which $f_\lambda$ is smooth. Then we can write the
sequence  of inequalities
\begin{eqnarray*}
\| f_\lambda\|_{ \dl^2_{\e}} &\lesssim& (\e \lambda)^{-1} \e^\frac{n}2
\| f_\lambda\|_{ \dl^\infty_{\lambda^{-1}}}
\\
 &\lesssim& (\e \lambda)^{-1} \e^\frac{n}2 \lambda^{\frac{n}2}
\| f_\lambda\|_{ \dl^2_{\lambda^{-1}}}
\\
 &\lesssim& (\e \lambda)^{-1} \e^\frac{n}2 \lambda^{\frac{n}2} \lambda
\| f_\lambda\|_{ \dl^2_{\lambda}}
\\
&=& \e^{\frac{n}2-1} \lambda^{\frac{n}2}\| f_\lambda\|_{ \dl^2_{\lambda}}
\end{eqnarray*}
All these bounds are obtained by comparing tubes with the same
orientation $\omega$.  The first step uses Holder's inequality to
switch from the $\e$ scale to the $\lambda^{-1}$ scale; the
$\dl^\infty$ norm is defined in the same way as the $\dl^2$ norm.  The
second takes advantage of the localization at frequency $\lambda$.
Finally the third step uses Holder's inequality to switch from the
$\lambda^{-1}$ scale to the $\lambda$ scale.
Summing up with respect to $\lambda$ we obtain \eqref{7a},
\[
\| f_{<\e^{-1}}\|_{\dl^2_{\e}} \lesssim \!\! \sum_{1 \leq \lambda \leq
  \e^{-1}}  \| f_\lambda\|_{{\dl^2_{\e}}} \lesssim 
\e^{\frac{n}2-1} \sum_{1 \leq \lambda \leq
  \e^{-1}} \lambda^{\frac{n}2}\| f_\lambda\|_{ \dl^2_{\lambda}} \lesssim \e^{\frac{n}2-1}\|f\|_{\mathcal{D} H^s}
\]
for $s > n/2$. 

Consider now the second part of \eqref{ulow}. The analogue of
\eqref{dscale} is
\[
\| \nabla f^\epsilon\|_{\dl^2_\lambda} = \e^{2-\frac{n}2} \| \nabla f\|_{\dl^2_{\e\lambda}} 
\]
therefore we have to show that
\[
\e^{2-\frac{n}2}\| \nabla f_{<\e^{-1}}\|_{\dl^2_{\e}} \lesssim \max\{\e,\e^{s-\frac{n}2}\}
\|f\|_{\mathcal{D} H^s}
\]
We proceed as above,
\[
\| \nabla f_{<\e^{-1}}\|_{\dl^2_{\e}} \lesssim \sum_{1 \leq \lambda \leq
  \e^{-1}}  \lambda \| f_\lambda\|_{{\dl^2_{\e}}} \lesssim 
\e^{\frac{n}2-1} \sum_{1 \leq \lambda \leq
  \e^{-1}} \lambda^{\frac{n}2+1}\| f_\lambda\|_{ \dl^2_{\lambda}}
\]
The bound for the last sum is straightforward and depends on the relative positions of $s$
and $\frac{n}2+1$.

It remains to consider the high frequency  estimate \eqref{uhigh}, for
which we only need the scaling relation \eqref{dscale} and Holder's inequality:
\[
\begin{split}
\| f^\e_{ > 1}\|_{\dh^s}^2 &= \sum_{\lambda > 1}  \lambda^{2s} \| (f^\e)_{\lambda}\|_{\dl^2_\lambda}^2
=  \e^{2-n} \sum_{\lambda > 1}  \lambda^{2s} \| f_{\e^{-1}
  \lambda}\|_{\dl^2_{\e \lambda}}^2
\\ & \lesssim \e^{2-n} \sum_{\lambda > 1} \e^{-2}  \lambda^{2s} \| f_{\e^{-1}
  \lambda}\|_{\dl^2_{\e^{-1} \lambda}}^2 = \e^{2s-n}  (\e^{-1} \lambda)^{2s}
\\ &
\lesssim \e^{2s-n} \|f\|_{\dh^s}^2
\end{split}
\]

\end{proof}

\begin{proof}[Proof of Proposition \ref{low2}] It suffices to consider
the $\tilde{a}(t,x,\xi)= u_{0, \leq 1}^{\e} \xi$ part of the symbol,
since $ \nabla{u_{0, \leq 1}^{\e}}$ is obtained from $\tilde{a}$
by differentiation. The main ingredient of the proof is

\begin{l1} \label{la} Let $s > \frac{n}2$. Then the following estimate
  holds:
  \beq \label{i1} 
\sup_{x_0 \in \R^n} \sup_{|\o_{0}|=1}  \int_{-\infty}^{\infty} |f(x_{0}+t
  \o_{0})| d t \leq \|f\|_{\dh^{s}} 
\eeq
\end{l1}

\begin{proof}[Proof of Lemma \ref{la}]

For each dyadic component $f_\lambda$ of $f$
we have:
\[
\begin{split}
  \int_{-\infty}^{\infty} |f_\lambda(x_0 + t \omega_0)| dt & = 
\sum_{k \in \Z} \int_{\lambda k }^{\lambda (k+1)} |f_\lambda(x_0 + t \omega_0)| dt \\
  &\les \lambda^\frac12 \sum_{k \in \Z} \left(
    \int_{\lambda k}^{\lambda (k+1)} 
|f^{2}_\lambda (x_0 + t \omega_0)| dt \right)^{\q} \\
  &\les \lambda^{\frac{n}{2}} \sum_{k \in \Z} \
\| 1_{\{|\lambda^{-1} (x-x_0)-k\omega|<1\}} f_\lambda \|_{L^2} \\
  &\les \lambda^{\frac{n}{2}} \|f_\lambda\|_{\dl^2_\lambda}
\end{split}
\]
At the second step we have used Holder's inequality and at the third
we have used the frequency localization.  The summation with respect
to $\lambda$ gives the desired result in \eqref{i1}.

 \end{proof}

 We return to the proof of the proposition.  A change of variables
 combined with Lemma~\ref{la} and the first part of \eqref{ulow} gives
\[
\sup_{x,\xi}\! \int_{0}^1 \! |\tilde a(t,x+2t\xi,\xi)| dt=\sup_{x,\xi}\!
\int_{0}^{2|\xi|}\! |u_{0, \leq 1}^{\e}(x+t\o_0)| dt \les \! \|u_{0, \leq
  1}^{\e} \|_{\mathcal{D}H^{s}} \les M
\]
where $\o_0=\frac{\xi}{|\xi|}$. 

For $x$ derivatives of $\tilde a$ we use the second part of
\eqref{ulow} instead:
\[
\int_0^{1} |\partial_x  \tilde{a} (x+2t \xi, \xi)| dt \les \| \nabla{u_{0, \leq 1}^{\e}}\|_{\mathcal{D}H^s} \les \max\{\e, \e^{s-\frac{n}{2}}\} M 
\]

For $\xi$ derivatives we use the pointwise bound \eqref{i2}:
\[
\int_0^{1} |\partial_\xi  \tilde{a} (x+2t \xi, \xi)| dt = \int_0^{1} |
u^\e_{0,\leq 1}(x+2t \xi)| dt \lesssim \e M
\]
 Since $u_{0, \leq 1}^{\e}$ is supported at frequencies $\les 1$, 
we obtain similar bounds for higher order derivatives.
\end{proof}

\begin{proof}[Proof of Proposition \ref{lowlin}] 
We neglect the gradients applied to $ u_{0, \leq 1}^{\e}$. Also for
the estimate involving $\bar v$ we retain only the stronger
bound\footnote{This is because we do not differentiate between
frequencies $1$ and less than $1$ in what follows}
involving $\nabla \bar v$. Then we need to show that
\begin{equation}
\| u_{0, \leq 1}^{\e} v\|_{Y^s} \lesssim \| u_{0, \leq
  1}^{\e}\|_{L^\infty}\| v\|_{X^s}
\label{uev}\end{equation}
respectively
\begin{equation}
 \| u_{0, \leq 1}^{\e}
\nabla \bar v\|_{Y^s} \lesssim \| u_{0, \leq
  1}^{\e}\|_{L^\infty}\| v\|_{X^s}.
\label{uevb}\end{equation}

In \eqref{uev} we use orthogonality with respect to unit frequency
cubes to reduce it to
\[
\| u_{0, \leq 1}^{\e} v_{\xi_0}\|_{Y} \lesssim \| u_{0, \leq
  1}^{\e}\|_{L^\infty}\| v_{\xi_0}\|_{X}
\]
By duality this becomes
\[
\left|\int  u_{0, \leq 1}^{\e} v_{\xi_0}  \bar w_{\xi_0} dx dt \right|
\lesssim\| u_{0, \leq
  1}^{\e}\|_{L^\infty}  \|v_{\xi_0}\|_{X} \|w_{\xi_0}\|_{X} 
\]
Hence it suffices to show that
\[
\int |  v_{\xi_0}|  | w_{\xi_0} |  dx dt  \lesssim
\|v_{\xi_0}\|_{X} \|w_{\xi_0}\|_{X} 
\]
By orthogonality with respect to $T_{\xi_0}$ tubes (see
Proposition~\ref{rep}) it remains to verify that given a $\xi$ tube $Q$
we can integrate a bump function on $Q$,
\[
\int 1_Q dx dt \lesssim 1 
\]
which is trivial.

Similarly, \eqref{uevb} reduces to  the frequency localized bound
\[
\left|\int  u_{0, \leq 1}^{\e} v_{\xi_0}  w_{- \xi_0} dx dt \right|
\lesssim |\xi_0|^{-1} \| u_{0, \leq
  1}^{\e}\|_{L^\infty}  \|v_{\xi_0}\|_{Y} \|w_{\xi_0}\|_{Y} 
\]
We use a modulation decomposition of $ v_{\xi_0}$ and   $w_{- \xi_0}$
at modulation $\xi^2/8$. Due to the Fourier localization, the 
integral corresponding to the low modulation parts vanishes,
\[
\int  u_{0, \leq 1}^{\e} M_{<\xi_0^2/4} v_{\xi_0}    M_{<\xi_0^2/4}
w_{- \xi_0} dx dt =0
\]
On the other hand we use Proposition~\ref{mod} to estimate
\[
\begin{split}
\left|\int  u_{0, \leq 1}^{\e} M_{>\xi_0^2/4} 
v_{\xi_0}   w_{-\xi_0} dx dt \right|
\lesssim & \|u_{0, \leq 1}^{\e}\|_{L^\infty} \|  M_{>\xi_0^2/4} 
v_{\xi_0}\|_{L^2} \| \|w_{-\xi_0}\|_{L^2}
\\ \lesssim &  |\xi_0|^{-1} \|u_{0, \leq 1}^{\e}\|_{L^\infty} \| 
v_{\xi_0}\|_{X^0}  \|w_{-\xi_0}\|_{X^0}
\end{split}
\]

\end{proof} 

\begin{proof}[Proof of Proposition \ref{low3}] This is a direct
  consequence of the estimates in \eqref{ulow} and \eqref{i2}.
\end{proof}

\section{The bilinear estimate}

Here we prove Theorem~\ref{ta3}.  After a Littlewood-Paley
decomposition it suffices to prove a bound for the high-low frequency
interactions
\begin{equation}
  \|  u_\mu  v_\lambda \|_{\D Y_\lambda} \lesssim \lambda^{-1}
  \mu^{\frac{n}2} (\ln \mu)^\frac12 \| u_\mu \|_{\D
    X_\mu}    \| v_\lambda \|_{\D X_\lambda}  \qquad 1 \leq \mu \ll \lambda
  \label{hlh}\end{equation}
respectively\footnote{Strictly speaking we should consider products of
  the form $u_{\lambda_1} v_{\lambda_2}$ with $\lambda_1 \approx
  \lambda_2$ but this makes no difference} for high-high frequency
interactions
\begin{equation}
  \| S_\mu (u_\lambda v_\lambda) \|_{\D Y_\mu} \lesssim \lambda^{n-1}
  \mu^{- \frac{n}2} \| u_\lambda \|_{\D
    X_\lambda}    \| v_\lambda \|_{\D
    X_\lambda}  \qquad 1 \leq \mu \lesssim \lambda
  \label{hhl}\end{equation}
and the similar bounds where one or both of the factors are replaced
by their complex conjugates.

Both bounds can be localized on the $2^j$ spatial scale. Thus the two
$D_j$'s can be factored out in the first bound, and neglected in the
second one. Then, by the duality result in Proposition~\ref{pdual},
\eqref{hlh} can be rewritten as
\begin{equation}
  \left| \int  u_\mu v_\lambda  \overline{w_\lambda} dx dt\right| \lesssim
  \lambda^{-1} \mu^{\frac{n}2} (\ln \mu)^\frac12  \|u_\mu \|_{\D
    X_\mu}    \| v_\lambda \|_{ X_\lambda}  \| w_\lambda \|_{ X_\lambda}
  \label{dual}\end{equation}

On the other hand, in \eqref{hhl} we can replace the $l^1$ summation
by an $l^2$ summation on the $\mu$ spatial scale by losing a
$\lambda^\frac12 \mu^{-\frac12}$ factor. Thus it remains to show that
\[
\| S_\mu (u_\lambda v_\lambda) \|_{Y_\mu} \lesssim \lambda^{n-\frac32}
\mu^{-\frac{(n-1)}2} \| u_\lambda \|_{X_\lambda} \| v_\lambda \|_{
  X_\lambda}
\]
By the duality result in Proposition~\ref{pdual} this is equivalent to
\[
\left| \int u_\lambda v_\lambda \overline{w_\mu} dx dt\right| \lesssim
\lambda^{n-\frac32} \mu^{-\frac{(n-1)}2} \| u_\lambda \|_{X_\lambda}
\| v_\lambda \|_{ X_\lambda}\|w_\mu \|_{X_\mu}
\]
which is easily seen to be weaker than \eqref{dual}.

It remains to prove \eqref{dual} where we allow $1 \leq \mu \leq
\lambda$ in order to include both cases above, and where we allow any
combination of complex conjugates.  The seemingly large number of
cases is reduced by observing that the bound rests unchanged if we
conjugate the entire product. Hence we can assume without any
restriction in generality that at most one factor is conjugated. Hence
it suffices to consider the following three cases:

(i) The product $u_\mu v_\lambda \overline{w_\lambda}$ with $1 \leq
\mu \leq \lambda$. This is the main case, where all three factors can
simultaneously concentrate in frequency near the parabola.

(ii) The product $\overline{u_\mu} v_\lambda {w_\lambda}$ with $1 \leq
\mu \ll \lambda$. Because the high frequency factors cannot
simultaneously concentrate on the parabola, the estimate turns
essentially into a bilinear $L^2$ estimate.

(iii) The product $u_\mu v_\lambda {w_\lambda}$ with $1 \leq \mu
\lesssim \lambda$. This is very similar to the second case.

{\bf Case 1:} Here we prove \eqref{dual} exactly as stated, for $1
\leq \mu \leq \lambda$.  This follows by summation with respect to
$\xi$, $\eta$ in the following result:

\begin{p1}
  a) Let $1 \leq \mu \leq \lambda$ and $\xi,\eta \in \R^n$ with
  \[
  |\xi|\approx \lambda, \qquad |\eta| \approx \mu, \qquad |\xi+\eta|
  \approx \lambda.
  \]
  Then we have
  \begin{equation}
    \begin{split}
      \left| \int S_\eta u S_\xi v \overline{S_{\xi+\eta} w} dx
        dt\right| \lesssim \frac{\| S_\eta u\|_{ X_{\eta}} \|S_\xi
        v\|_{X_\xi} \| S_{\xi+\eta} w\|_{X_{\xi+\eta}}}{
        \mu^{\frac{1}2} (\lambda + \min\{|\xi\cdot \eta|, |\xi \wedge
        \eta|\})^{\frac12}}
    \end{split}
    \label{las}\end{equation}

  b) Assume in addition that $\mu \lesssim \lambda^\frac12$. Then we
  have
  \begin{equation}
    \begin{split}
      \left| \int S_\eta u S_\xi v \overline{S_{\xi+\eta} w} dx
        dt\right| \lesssim \frac{\| S_\eta u\|_{\D_i X_{\eta}} \|S_\xi
        v\|_{X_\xi} \| S_{\xi+\eta} w\|_{X_{\xi+\eta}}}{
        \lambda^{\frac{1}2} \mu^{-\frac{1}2} (\lambda +|\xi\cdot
        \eta|)^{\frac12} }
    \end{split}
    \label{las1}\end{equation}
  \label{ptubes}\end{p1}

We first show how to use the proposition to conclude the proof of
\eqref{dual}. We decompose $u_\mu$, $v_\lambda$ and $w_\lambda$ in
unit frequency cubes,
\[
u_\mu = \sum_{|\eta| \approx \mu} S_\eta u_\mu, \qquad v_\lambda =
\sum_{|\xi| \approx \lambda} S_\xi v_\lambda, \qquad w_\lambda =
\sum_{|\zeta| \approx \lambda} S_\zeta w_\lambda
\]
and use the corresponding decomposition of the integral in
\eqref{dual}.

We consider two cases. If $\mu \geq \lambda^\frac12$ then we use
\eqref{las} to estimate
\[
\left| \int u_\mu v_\lambda \overline{w_\lambda} dx dt\right| \lesssim
\sum_{|\xi| \approx \lambda,\ |\eta| \approx \mu}^{|\xi+\eta| \approx
  \lambda} \frac{\| S_\eta u_\mu \|_{ X_{\eta}} \|S_\xi
  v_\lambda\|_{X_\xi} \|S_{\xi+\eta} w_\lambda\|_{X_{\xi+\eta}}}{
  \mu^{\frac{1}2} (\lambda + \min\{|\xi\cdot \eta|, |\xi \wedge
  \eta|\})^{\frac12}}
\]
By the Cauchy-Schwartz inequality this is bounded by
\[
\left(\sum_{|\xi| \approx \lambda,\, |\eta| \approx \mu}^{|\xi+\eta|
    \approx \lambda} \| S_\eta u_\mu \|_{ X_{\eta}}^2 \|S_{\xi+\eta}
  w_\lambda\|_{X_{\xi+\eta}}^2\!\! \right)^\frac12 \!\!\!
\left(\sum_{|\xi| \approx \lambda,\, |\eta| \approx \mu}^{|\xi+\eta|
    \approx \lambda}\frac{\|S_\xi v_\lambda\|_{X_\xi}^2 } { \mu
    (\lambda + \min\{|\xi\cdot \eta|, |\xi \wedge \eta|\})} \!\!
\right)^\frac12
\]
and further by
\[
\|u_\mu \|_{ X_\mu} \| v_\lambda \|_{ X_\lambda} \| w_\lambda \|_{
  X_\lambda} \left( \sup_{|\xi| \approx \lambda} \sum_{|\eta| \approx
    \mu}^{|\xi+\eta| \approx \lambda} \frac{1} { \mu (\lambda +
    \min\{|\xi\cdot \eta|, |\xi \wedge \eta|\})}\right)^\frac12
\]
which gives \eqref{dual} since
\[
\sum_{|\eta| \approx \mu}^{|\xi+\eta| \approx \lambda} \frac{1} { \mu
  (\lambda + \min\{|\xi\cdot \eta|, |\xi \wedge \eta|\})} \approx
\frac{\mu^{n} \ln \mu}{\mu^2 \lambda}
\]
and $\mu^2 \geq \lambda$. We note that the bound improves as $\mu$
increases.

If $\mu \leq \lambda^\frac12$ then the argument is similar but using
\eqref{las1} instead of \eqref{las}. This concludes the proof of
\eqref{dual}.

\begin{proof}[Proof of Proposition~\ref{ptubes}]
  We decompose each of the factors in wave packets,
  \[
  S_\eta u = \sum_{P \in T_\eta} u_P, \qquad S_\xi v = \sum_{Q \in
    T_\xi} v_Q, \qquad S_{\xi+\eta} w = \sum_{R \in T_{\xi+\eta}} w_R
  \]
  We first prove a bound with $Q,R$ fixed.

\begin{l1}
  For $\xi$ and $\eta$ as above the following estimate holds:
  \begin{equation}
    \left| \int \sum_{P \in T_\eta}  u_P v_Q \overline{w_R} dx\right| \lesssim
    \left(\sum_{P \in T_\eta}^{P\cap Q \cap R \not= \emptyset}\|
      u_P\|^2_{X_\eta}\right)^\frac12 \frac{ \| v_Q\|_{X_\xi} \|w_R\|_{X_{\xi+\eta}}}{\mu^{\frac{1}2} 
      (\lambda +|\xi\cdot \eta|)^{\frac12} }
  \end{equation}
  \label{nodl}\end{l1}

\begin{proof}
  By Proposition~\ref{rep} we can assume without any restriction in
  generality that
  \[
  u_P = a_P(t) e^{-it\eta^2} \chi_P(x-2t\eta), \qquad v_Q = a_Q(t)
  e^{-it\xi^2}\chi_Q(x-2t\xi),
  \]
  \[
  \overline{w_R} = a_R(t) e^{it(\xi+\eta)^2} \chi_R(x-2t(\xi+\eta))
  \]
  with the $\chi$'s being unit bumps and the $a$'s in $V^2$.  Then the
  integral has the form
  \[
  \int_{-1}^1 \!\!\!  e^{4it \xi \eta} a_Q(t) a_R(t) \! \sum_P \! a_P(t)
  \!\! \int_{\R^n} \!\!\!\!  \chi_P(x-2t\eta) \chi_Q(x-2t\xi)
  \chi_R(x-2t(\xi+\eta)) dx dt
  \]

  The tubes $Q$ and $R$ differ in speed by $\eta$, therefore they
  intersect in a time interval $I$ of lenght at most $\mu^{-1}$.  The
  tubes $P$ and $R$ differ in speed by $\xi$, therefore they intersect
  in a time interval of lenght at most $\lambda^{-1}$.  Thus there are
  about $\lambda \mu^{-1}$ tubes $P$ which intersect both $Q$ and $R$.

  For fixed $P$ the $x$ integral above is a smooth bump function on a
  $\lambda^{-1}$ interval. Then we can express the above integral in
  the form
  \[
  \int_I e^{4it \xi \eta} a_Q(t) a_R(t) \sum_P a_P(t) b_P(t) dt
  \]
  where $b_P$ are smooth bump functions on essentially disjoint
  $\lambda^{-1}$ intervals. For the sum with respect to $P$ we can
  estimate
  \[
  \| \sum _{P} a_P(t) b_P(t)\|_{L^2}^2 \lesssim \lambda^{-1} \sum \|
  a_P\|_{L^\infty}^2
  \]
  and
  \[
  \| \sum _{P} a_P(t) b_P(t)\|_{V^2}^2 \lesssim \sum \| a_P\|_{V^2}^2
  \]
  We consider two possibilities. If $|\xi\cdot \eta| \lesssim \lambda$
  then we use Holder's inequality to bound the integral by
  \[
  \mu^{-\frac{1}2} \lambda^{-\frac12} \sum_P \| a_P\|_{L^\infty}
  \|a_Q\|_{L^\infty} \|a_R\|_{L^\infty}
  \]
  If $|\xi \cdot \eta| \gtrsim \lambda$ then we use the algebra
  property for $V^2$. It remains to prove that
  \[
  \left |\int_{0}^{\mu^{-1}} a(t) e^{i t \sigma} dt \right| \lesssim
  \mu^{-\frac{1}2} |\sigma|^{-\frac12} \|a\|_{V^2} \qquad \sigma = 4
  \xi \eta
  \]
  After rescaling this becomes
  \[
  \left| \int_{0}^1 a(t) e^{i t \sigma} dt \right| \lesssim
  |\sigma|^{-\frac12} \|a\|_{V^2}
  \]
  This follows by Holder's inequality from the bound
  \[
  \| S_{\geq \sigma} a(t)\|_{L^2} \lesssim \sigma^{-\frac12}
  \|a\|_{V^2}
  \]

\end{proof}

Now we prove part (a) of the proposition.  Three tubes $P,Q,R$
contribute to the integral only if they intersect.  We consider the
intersection pattern of $\xi$ and $\xi+\eta$ tubes.  For any $\xi$
tube $Q$, all $\xi+\eta$ tubes intersecting it are contained in a
larger slab obtained by horizontally translating $Q$ in the $\eta$
direction
\[
H = 2(Q + \{0\} \times [-2\eta,2\eta])
\]
We denote by $\H$ a locally finite covering of $[-1,1] \times \R^n$
with such slabs. Then we use the lemma to bound the integral in
\eqref{las} by
\[
\mu^{-\frac{1}2} (\lambda+|\xi\cdot \eta|)^{-\frac12} \sum_{H \in \H}
\sum_{Q,R \subset H} \|u_Q\|_{X_\xi} \| u_R\|_{X_{\xi+\eta}}
\left(\sum_{P \in T_\eta}^{P\cap Q \cap R \not= \emptyset}\|
  u_P\|^2_{X_\eta}\right)^\frac12
\]
Using Cauchy-Schwartz in the second sum with respect to $(Q,R)$ we
bound this by
\[
\frac{N^\frac12}{\mu^{\frac{1}2} (\lambda+|\xi\cdot \eta|)^{\frac12}}
\left(\sum_{P \in T_\eta}\| u_P\|^2_{X_\eta}\right)^\frac12 \sum_{H
  \in \H} \left(\sum_{Q \subset H} \|u_Q\|_{X_\xi}^2\right)^\frac12
\left(\sum_{R \subset H} \| u_R\|_{X_{\xi+\eta}}^2\right)^\frac12
\]
where
\[
N = \max_{H,P} |\{ (Q,R);\ Q,R \subset H,\ P\cap Q \cap R \not=
\emptyset\}|
\]
An additional Cauchy-Schwartz allows us to estimate the above sum by
\[
\frac{N^\frac12}{ \mu^{\frac{1}2} (\lambda+|\xi\cdot \eta|)^{\frac12}}
\left(\sum_{P \in T_\eta}\| u_P\|^2_{X_\eta}\right)^\frac12
\left(\sum_{Q } \|u_Q\|_{X_\xi}^2\right)^\frac12 \left(\sum_{R} \|
  u_R\|_{X_{\xi+\eta}}^2\right)^\frac12
\]
To conclude the proof it remains to establish a bound for $N$, namely
\[
N \lesssim \frac{|\xi||\eta|}{|\xi| + |\xi \wedge \eta|}
\]
Both $P \cap R$ and $P \cap Q$ have a $\lambda^{-1}$ time length and
are uniquely determined by this intersection up to finite
multiplicity.  It follows that
\[
N \lesssim \lambda |I|
\]
where $I$ is the time interval where $H$ and $P$ intersect.  Given the
definition of $H$ it follows that $I$ has the form
\[
I = \{ t ; |x_0 + t \xi + s \eta| \leq 2 \ \text { for some } s \in
[-1,1]\}
\]
Symmetrizing we can assume that $x_0=0$. Taking inner and wedge
products with $\eta$ it follows that $t$ must satisfy
\[
|t| \leq |\eta|^2 |\xi \cdot \eta|^{-1}, \qquad |t| \leq |\eta| |\xi
\wedge \eta|^{-1}
\]
which lead to the desired bound for $N$.

b) We first note that both $Q$ and $R$ are contained in spatial strips
of size $\mu \times \lambda$ oriented in the $\xi$ direction. Hence by
orthogonality it suffices to prove the estimate in a single such
strip. Then we can take advantage of the $l^1$ summability in the $\D
X_\mu$ norm to further reduce the estimate to the case when the $\eta$
tubes are spatially concentrated in a single $\mu \times \mu$ cube
$Z$.

A $\xi$ or a $\xi+\eta$ tube needs a time of $\mu \lambda^{-1}$ to
move through such a cube. On the other hand, the tubes $Q$ and $R$
need a larger time $\mu^{-1}$ to separate. Hence within $Z$ we can
identify the $\xi$ tubes and the $\xi+\eta$ tubes. By orthogonality it
suffices to consider a single $Q$ and a single $R$.  Consequently, the
conclusion follows from the following counterpart of Lemma~\ref{nodl}.

\begin{l1}
  For $\xi$ and $\eta$ as in part (b) of the proposition the following
  estimate holds
  \begin{equation}
    \left| \int \sum_{P \in T_\eta}^{P \subset Z}  u_P v_Q w_R dx\right| \lesssim
    \left(\sum_{P \in T_\eta}^{P\cap Q \cap R \not= \emptyset}\|
      u_P\|^2_{X_\eta}\right)^\frac12 \frac{\| v_Q\|_{X_\xi}
      \|w_R\|_{X_{\xi+\eta}}}{
      \lambda^{\frac{1}2} \mu^{-\frac{1}2} (\lambda +|\xi\cdot \eta|)^{\frac12}}  
  \end{equation}
  \label{dl}\end{l1}

The proof is almost identical with the proof of Lemma~\ref{nodl}, the
only difference is that we now have $|I| \approx \mu \lambda^{-1}$.

{\bf Case 2:} Here we consider the product $\overline{u_\mu} v_\lambda
{w_\lambda}$ where $1 \leq \mu \ll \lambda$ and prove a stronger
bound, namely
\begin{equation}
  \left| \int   {\overline u_\mu} v_\lambda {w_\lambda} dx dt\right| \lesssim
  \lambda^{-1} \mu^{\frac{n}2-1}   \|u_\mu \|_{
    X_\mu}    \| v_\lambda \|_{ X_\lambda}  \| w_\lambda \|_{ X_\lambda}
  \label{dualgood}\end{equation}

We use the modulation localization operators $S_{< \lambda^{2} /100}$ to
split each of the factors in two,
\[
u_\mu = M_{< \lambda^{2} /100} u_\mu + (1-M_{< \lambda^{2} /100} u_\mu)
\]
etc. We observe that
\[
\int \overline{ M_{< \lambda^{2} /100} u_\mu} M_{< \lambda^{2} /100}
v_\lambda { M_{< \lambda^{2} /100} w_\lambda} dx dt= 0
\]
due to the time frequency localizations. Precisely, the first factor
is frequency localized in the region $\{ |\tau| < \lambda^2/50\}$
while the other two are frequency localized in the region $\{ \tau >
\lambda^2/8\}$. Hence it remains to consider the case when at least
one factor has high modulation. For that factor we have a favourable
$L^2$ bound as in Proposition~\ref{mod},
\[
\| M_{> \lambda^{2} /100} u_\mu \|_{L^2} \lesssim \lambda^{-1} \|
u_\mu\|_{X_\mu}
\]
and similarly for the other factors. Then it remains to prove that
\[
\| \bar u_\mu v_\lambda\|_{L^2} \lesssim \mu^{\frac{n}2-1} \|u_\mu \|_{\D
  X_\mu} \| v_\lambda \|_{ X_\lambda}
\]
respectively
\[
\| S_\mu (v_\lambda {w_\lambda})\|_{L^2} \lesssim \mu^{\frac{n}2-1} \|
v_\lambda \|_{ X_\lambda} \| w_\lambda \|_{ X_\lambda}
\]
For this we need the Strichartz estimates in Proposition~\ref{SE} (see
also \eqref{set}).

For the first bound we use the $L^\frac{2(n+2)}n$ estimate for
$v_\lambda$, respectively the $L^{n+2}$ estimate for $\bar u_\mu$.

For the second we first oberve that by orthogonality it suffices to
prove it when both $v_\lambda$ and $w_\lambda$ are frequency localized
to cubes of size $\mu$.  Then we use the $L^\frac{2(n+2)}n$ estimate
for both factors to derive $L^4$ bounds by Sobolev embeddings.

{\bf Case 3:} Here we consider the product ${u_\mu} v_\lambda
{w_\lambda}$ where $1 \leq \mu \lesssim \lambda$. This is treated
exactly as above, and an estimate similar to \eqref{dualgood} is
obtained.
 
\end{proof}

\end{document}